\documentclass[journal]{IEEEtran}
\usepackage[noadjust]{cite}

\ifCLASSINFOpdf
\else
\fi
\usepackage{amsmath}
\usepackage{mathrsfs}
\usepackage{bm}
\usepackage{amssymb}

\usepackage{graphics}
\usepackage{makecell}
\usepackage{graphicx}
\usepackage{enumerate}

\begin{document}


\title{Controlling Time-Variant Virtual Inertia from Storage by Dynamic Programming and PROPT\\}

%\author{Shuchang~Yan,~\IEEEmembership{Student Member,~IEEE,}
%	Yu~Zheng,~\IEEEmembership{Member,~IEEE,}
%	and~David~John~Hill,~\IEEEmembership{Life~Fellow,~IEEE}% <-this % stops a space
%	\thanks{M. Shell was with the Department
%		of Electrical and Computer Engineering, Georgia Institute of Technology, Atlanta,
%		GA, 30332 USA e-mail: (see http://www.michaelshell.org/contact.html).}% <-this % stops a space
%	\thanks{J. Doe and J. Doe are with Anonymous University.}% <-this % stops a space
%	\thanks{Manuscript received April 19, 2005; revised August 26, 2015.}}

\author{Shuchang~Yan,~\IEEEmembership{Student Member,~IEEE}
	%Yu~Zheng,~\IEEEmembership{Member,~IEEE,}
	%and~David~John~Hill,~\IEEEmembership{Life~Fellow,~IEEE}% <-this % stops a space
	%\thanks{The work is supported by the Hong Kong RGC Theme Based Research Scheme under Grants T23-407/12N and T23-701/14N.}}
}

\maketitle

\begin{abstract}

The integration of renewables gradually replaces the traditional power plants, and this makes that the rotational inertia provided by the power plants is decreasing with time. Virtual inertia emulated by power electronic devices is becoming a promising way to stabilize the frequency of power system when a disturbance happens. And how to control virtual inertia to achieve desired performance is an important question to be answered. In this work, authors formulate an optimal control problem to describe the behaviour of controlling the time-variant inertia from storage, and structure preserving model is utilized to describe the dynamics of system, where frequencies of all buses are preserved. Also, practical constraints such as frequency contraints for buses, power/energy constraints for storage are incorporated in the optimal control problem. To solve this optimal control problem, dynamic programming (DP)  and PROPT MATLAB Optimal Control Software are employed to obtain global optimal solution (virtual inertia trajectory). The simulation results show the correctness and effectiveness of our problem formulation and solving method.

\end{abstract}

\IEEEpeerreviewmaketitle

\section{Introduction}

%Due to the economic and enviromental concern, many renewable energy sources such as wind turbines and photovoltaic cells have been introduced into the power system and high proportion of conventional generators based on fossil fuels have been replaced by renewable energy sources. The wind power has the characteristics of variability and uncertainty, and thus, the ramping of the net load is much higher than before, which is difficult for conventional generators to catch up. On the other hand, inertia and damping provided by conventional synchronous generators are decreasing, whether traditional frequency control strategies can ensure frequency stability of power system becomes a question.
%Current research considers the possibilities of utilizing other devices such as energy storage to mimic the behaviors of synchronous generators, that is, virtual inertia and damping. And the research question is {\itshape{where and how to control these devices can have the strongest support to the frequency of the power system}.} The corresponding research methods include the simulation-based method and model-based method. On the other hand, renewables connects to the power grid with power electronic devices, and most of these devices in practice do not provide inertia support

To achieve better enviromental and economic benefits, more and more renewable energy sources are brought into the current power system.  On the one hand, the wind power is uncertian and variable, and thus, wind power prediction with 100 percent accuracy is impossible. This will make power imbalance happens more frequently than before. On the other hand,  the traditional power plants are decreasing while most of renewables connect to the power grid with power electronic devices where no inertia is provided, which means a low-inertia power system is coming now~\cite{ulbig2014impact}. Researchers in power system area proposes to use power electronic devices to emulate virtual inertia to stabilize the frequency when a disturbance happens. And the research question is {\itshape{where to locate and how to control virtual inertia from these power electronic-interfaced devices can have the strongest support to the frequency of the power system}.}

To anwser this question, two factors, the location to implement the control and the control method for virtual inertia should be seriously considered. Reference~\cite{adrees2016study} investigates the effect of deploying energy storage with droop control centrally and distributedlly on the primary frequency regulation of power system.  Reference~\cite{adrees2017impact} finds that deploying energy storage with droop contol can improve the transient performance indices such as frequency nadir by comparing power system with and without energy storage. Reference~\cite{adrees2018influence} further extends the work in reference~\cite{adrees2016study} and focuses on the relationship between the location of energy storage and primary frequency response in the low-inertia system. It is noted that references~\cite{adrees2016study,adrees2017impact, adrees2018influence} adopt a pure simulation approach to analyze the question. The benefit is that we can obtain a deterministic result for the optimal location and control parameters of storage for primay frequency response, however, the dificiencies are also obvious: the solution obtained in a certian grid may not be suitable for other grids and the control parameters during the transient process are fixed.

To design a more scientific method to analyze this question, reference~\cite{borsche2015effects} adopts a gradient-based method to analyze how the virtual  inertia and damping allocate can have desired transient performance. Reference~\cite{borsche2017placement} takes one step forward to include the transfer function of energy storage into that of power system, and analyzes that the relationship between frequency indices and virtual inertia/damping. However, the modelling in references~\cite{borsche2015effects,borsche2017placement} fails to analzye the effects of the location of energy storage on the frequency indices. References~\cite{poolla2018virtual, poolla2018placement} utilize $H_{2}$ norm to describe the transient performance. The allocation result of virtual inertia and damping can ensure the coherency of system, but it is noted that frequency contraints such as lower and upper limits can not be implemented in their problem formulations. Reference~\cite{guggilam2018optimizing} takes the assumption that frequency dynamics at all buses are the same in transient process, and then utilizes the second-order model to replace the original power system model to give the analytical relationshiop between frequency indices and virtual inertia and damping. However, the assumption is reasonable for a meshed network, and it is not in line with reality for certain types of grid structure such as chain-type system. Also, there is a common deficiency in references~\cite{adrees2016study}-\cite{guggilam2018optimizing} that control parameters or virtual inertia/damping during the transient process are fixed.

For the research on time-variant inertia, reference~\cite{fang2018inertia} reviews the control methods for the power electronic equipments to emulate virtual inertia. To the best of author's knowledge, there is only one reference~\cite{markovic2018lqr} to analyze how to control time-variant inertia to meet frequency requirements. However, there are several deficiencies: first, high-order power system dynamics is approximated by second-order system, which means the frequencies of different buses (or stucture information) are covered by one aggregated frequency; second, the LQR optimization technique is based on the linearized model, so the virtual inertia trajectory may not be optimal or appropriate for a large disturbance of the power system.

Based on the literature review above, authors here point out the design requirements for problem formulation and the corresponding control methods of time-variant virtual inertia from power electronic devices:  
   \begin{enumerate}[i.]
	\item The problem formualtion should allow that frequency of all buses can be presented or analyzed.
	\item The control method can better give a global optimal rather than local optimal trajectory of control inputs. 
	\item The control method (or solving algorithm) should be suitable to solve different system dynamics, and this is due to the fact the dynamics of power electronic devices varies from each other. 
	\item The control method can make frequency trajectory meet the specific requirements such as lower and upper frequency limits in the transient process and the steady frequency requirements after the disturbance.
	\end{enumerate}

To meet the above requirements, authors formulate an optimal control problem to obtain the optimal time-variant inertia trajectory and dynamic programming, as a widely applicable method, is employed here to solve this optimal control problem. For high dimensional optimal control problem, PROPT Matlab Optimal Control Software is utilized here to deal with this optimal control problem.

The remainder of this paper is structured as follows. In Section II, the structure preserving model and energy storage model are presented. Section III formulates the optimal  problem including the objectives and constraints. Section IV illustrates that how to utilize dynamic programming and PROPT to solve this optimal control problem. The case study is done in section V to verify the problem formulation and the corresponding solving method. The last section (Section VI) gives the final conclusion and future research directions.

\section{Modeling}

\subsection{Power System Dynamics}

The structure preserving model~\cite{bergen1981structure} is utilized here to model the dynamics of power system. In this model, reactive power is ignored and voltage magnititudes at all buses are assumed to be constant at 1 per unit value (p.u.). And all the transmission lines in the power system are assumed to be resistanceless. There are $N$ nodes or buses 	in the power system. Nodes or buses with inertia and without inertia are denoted as the set $G$ and set $L$ respectively, and superscripts $g$ and $l$ are the elements from $G$ and $L$ respectively. In power system, the nodes with inertia are generally buses with generators or motors, and the corresponding dynamics are described as follows, 
\begin{align}
\dot{\delta_{i}^{g}}=\omega_{i}^{g} \label{eq1}\\
\dot{\omega}_{i}^{g}=-\frac{D_{i}}{M_{i}}\omega_{i}^{g}-\frac{1}{M_{i}}&\sum^{N}_{\substack{j=1 \\ j  \neq i}}b_{ij}\sin(\delta_{i}-\delta_{j})+\frac{1}{M_{i}}P_{i}^{0} \label{eq2}
\end{align}
where $\delta_{i}$ and $\delta_{j}$ are the angle of bus $i$ and bus $j$ respectively, and $M_{i}$ and $D_{i}$ are the inertia and damping of bus $i$ respectively, $\omega_{i}$ is the frequency  at bus $i$, and it takes the nominal frequency as a reference, $b_{ij}$ is the susceptance between bus $i$ and bus $j$, and $b_{ij}\sin(\delta_{i}-\delta_{j})$ is the active power flow from bus $i$ to bus $j$, $P_{i}^{0}$ is a shorthand for $P_{M,i}^{0}-P_{D,i}^{0}$, which is the difference between mechanical power input $P_{M,i}^{0}$ of generator at bus $i$ and load demand $P_{D,i}^{0}$ at bus $i$. For the nodes without inertia, usually load buses without motor loads, the dynamics can be descibed as follows,
\begin{align}
\dot{\delta}_{i}^{l}=-\frac{1}{D_{i}}&\sum^{N}_{\substack{j=1 \\ j  \neq i}}b_{ij}\sin(\delta_{i}-\delta_{j})+\frac{1}{D_{i}}P_{i}^{0} \label{eq3}
\end{align}
where the mechanical power input $P_{M,i}^{0}$ at buses without inertia is equal to 0 and $P_{i}^{0}$ is equal to $-P_{D,i}^{0}$, representing the active power drawn from the node $i$.

\subsection{Energy Storage Dynamics}

Energy storage can mimic the behaviour of synchronous generators to provide virtual inertia and damping, the dynamics of energy storage is given below~\cite{guggilam2018optimizing},
\begin{align}
\dot{\delta_{i}^{e}}=\omega_{i}^{e} \label{eq4}\\
\dot{\omega}_{i}^{e}=-\frac{D_{e, i}}{M_{e,i}}\omega_{i}^{e}-\frac{1}{M_{e, i}}&\sum^{N}_{\substack{j=1 \\ j  \neq i}}b_{ij}\sin(\delta_{i}-\delta_{j})+\frac{1}{M_{e, i}}P_{i}^{e} \label{eq5}
%\dot{E}_{i}=P_{i}^{e}
\end{align}
where $M_{e, i}$ and $D_{e, i}$ are the virtual inertia and damping for the energy storage at bus $i$ respectively, $P_{i}^{e}$ is the constant power input or power output for energy storage at bus $i$. The nodes or buses with energy storage are denoted with the set $S$. For simplicity's consideration, there exists $N=G\cup L \cup S$, $ G \cap L=\varnothing$, $ G \cap S=\varnothing$ and $ L \cap S=\varnothing$.

\subsection{Constraints for the System Dynamics}

The dynamics of power system with energy storage can be expressed by~(\ref{eq1})-(\ref{eq5}). However, there are some practical limits such as the frequency limits for several certain buses. The system constraints are listed as follows,
\begin{gather}
|\omega_{i}| \leq \omega_{i}^{max}\label{eq6}\\
\omega_{i}(t_{1})~is~within~a~predefined~range\label{eq7}\\
P_{i}^{r,min}\leq P_{i}^{r}= P_{i}^{e}-M_{e,i}\dot{w^{e}}_{i}-D_{e,i}w^{e}_{i}\leq P_{i}^{r,max}\label{eq8}\\
E^{al,l}_{i}\leq \int_{t_{0}}^{t_{1}}P_{i}^{r} dt\leq E^{al,u}_{i}\label{eq9}
\end{gather}
where $\omega_{i}^{max}$ is the allowable frequency change of bus $i$, $P_{i}^{r,min}$ and $P_{i}^{r,max}$ are the lower and upper power limits for energy storage at bus $i$, and $E^{al,l}_{i}$ and $E^{al, u}_{i}$ are the lower and upper limits of energy change for storage at bus $i$, and $t_{1}$ is final time instant of concerned time interval. It should be noted that the specific values of $E^{al, l}$ and $E^{al, u}$ should be based on the current state of charge of storage.

\section{Problem Formulation}

%As is described in the introduction part, the aim of this work is to obtain the desired transient process including reducing the frequency overshoot,  energy exchange of storage and increasing the decay rate of oscillations by approriately locating and allocating virtual inertia and damping from storage. And thus, the optimization problem is formulated as follows,

The aim of this work is to control the virtual inertia to make sure that transient performance can be enhanced while the power and energy changes of storage do not go beyond their limits, and also specific frequency requirements of certain buses could be met during the transient process.

%\begin{gather}
%\label{eq24}
%\min_{\substack{{x_{q},m_{q},k_{q}}}}-c_{1}\xi^{min}+c_{2}\omega_{p}^{max}+c_{3}\dot{\omega}_{p}^{max}\nonumber\\+c_{4}\sum_{q\in Q}\Delta E_{e,q}+pen
%\end{gather}
\vskip 0.3cm

\begin{gather}
\min_{\substack{{M_{e,s}(t)}}}\int_{t_{0}}^{t_{1}}\sum_{{s}\in S}a_{s}(M_{e,s}(t)-M_{e,i,d})^{2}\nonumber \\+\sum_{{i}\in N}b_{i}|\omega_{i}(t)|\nonumber\\
+\sum_{{i}\in N}c_{i}\delta_{i}(t)^{2}~dt \label{eq10}
\end{gather}

\vskip 0.3cm
subject to
\begin{gather}
%q\in Q \label{eq25-}\\
System~dynamic~constraints~~~(\ref{eq1})-(\ref{eq5})\nonumber\\
~~~~~~~Frequency~contraints~~~(\ref{eq6}), (\ref{eq7})\nonumber\\
~~Power/energy~contraints~~~(\ref{eq8}), (\ref{eq9})\nonumber
%M_{e,i}(t_{1})=M_{e,i}(0)~~~~~~~~~~(9) \nonumber
\end{gather}
where $M_{e,s,d}$ is desired or reference virtual inertia value provided by storage at bus $s$; $M_{e,s}(t)$, virtual inertia provided by storage at bus $s$, is the control input; $\delta_{i}$ and $\omega_{i}$, angle and frequency respectively, are the state variables. We denote the $M_{e,s}(t)$ as $u_{e,s}(t)$, denote $\delta_{i}$ and $\omega_{i}$ as $x_{i}(t)$, and denote the term to be integrated in the objective function as $g(\mathbf{x(t)},\mathbf{u(t)})$ for convenience, where $\mathbf{x(t)}$ and $\mathbf{u(t)}$ are stacked by $x_{i}(t)$ and $u_{e,s}(t)$ respectively. And then, the optimal control problem could be expressed as the following concise form,

\begin{gather}
\min_{\substack{{\mathbf u(t)}}}\int_{t_{0}}^{t_{1}}g(\mathbf{x}(t),\mathbf{u}(t))dt\label{eq11}\\
\dot{\mathbf x}(t)=\mathbf f({\mathbf x}(t),{\mathbf u}(t))\label{eq12}\\
{\mathbf x}(t) \subset {\mathbf X}(t) \label{eq13}\\
{\mathbf u}(t) \subset {\mathbf U}(t) \label{eq14}\\
Other~~constraints \label{eq15}
\end{gather}
where ${\mathbf X}(t)$ is allowable trajectory for ${\mathbf x}(t)$, and ${\mathbf U}(t)$ is allowable trajectory for ${\mathbf u}(t)$.
The objective function is to minimize the virtual inertia change (control effort) and angle/frequency change (control performance) of all buses in the system. %Constraint (9) is to limit that the final value of virtual inertia after control implementation comes back to the initial value, which is practical and desired in the engineering practice. 
\section{Solving Method}
%The ridDE algorithm is adopted in our work and the benefit is that we can optimze the location, and the virtual inertia and damping from energy storage simultaneously. It is modified to solve the problem proposed in our work as follows,
In this section, dynamic programming and PROPT MATLAB OPTIMAL CONTROL SOFTWARE will be introduced respectively to solve this optimal control problem. 
\subsection{Dynamic Programming}
%\vskip 0.3cm

Dynamic programming is a classical method to solve optimal control problem, and its basic idea is to break a complex optimization problem into multi-stage subproblems, and the full-stage global optimal solution (or trajectory) is obtained through combining the optimal solution of subproblems, which is also known as ``Bellman's Principle of Optimality"~\cite{bellman1965dynamic, bertsekas1995dynamic}. The benefit is that we can obtain a global optimal trajectory within constraints for an optimal control problem. 

In this work, we adopt level-set dynamic programming (LS DP) proposed in reference~\cite{elbert2013implementation} to solve this problem and the results are also compared with these solved by basic DP. The difference of LS DP and the basic DP algorithm is that the former one utilizes interpolation between backward-reachable and non-backward-reachable grid points so that the obtained solution accuracy will be increased dramatically. The procedures to implement the LS DP and basic DP into our optimal control problem are listed as follows, and for the details of LS DP algorithm, readers please refer to reference~\cite{elbert2013implementation}:

\vskip 0.3cm

Step 1)\quad{\itshape{Initial Parameter Setting}}
\begin{itemize}
	\item[]
	\setlength\leftskip{4em}
	\quad Choose a node/bus for angle/frequency reference in the power system. This is quite important because the value of angle and frequency of other buses will not fly to the distance if we have a reference bus.
	
	\quad Choose the time interval for the simulation ($t_{0}$ and $t_{1}$ respectively), the lower and upper bounds of states variables ${\mathbf x}(t)$ and control inputs ${\mathbf u}(t)$, the value or ranges for the state variables ${\mathbf x}(t)$ at the final time $t_{1}$.
	
	\quad Choose the time step ($T_{s}$) for discretization of system dynamics, the penalty for final state ($\phi_{N}(\mathbf x_{N})$, also denoted as $My.Inf$ in Matlab code~\cite {kkk}) and number of state points in the grid ($N_{x}$). $N_{x}$ is utilized to discretize the space of state variables or control input $x$.
\end{itemize}

\vskip 0.1cm

Step 2)\quad{\itshape{Discretize the System Dynamics}}

\begin{itemize}
	\item[]
	\setlength\leftskip{4em}
	\quad Discretize the system dynamics based on the initial parameter setting as follows, 
	\begin{gather}
	\min_{\substack{{\mathbf u(k)}}}g_N(\mathbf x_{N})+\sum_{N_{0}}^{N-1}(g_{k}(\mathbf x_{k},\mathbf u_{k})\times T_{s})\label{eq16}\\	
	{\mathbf x}(k+1)=\mathbf f_{k}({\mathbf x}(k),{\mathbf u}(k))\label{eq17}\\
	{\mathbf x}(k) \subset {\mathbf X}(k) \label{eq18}\\
	{\mathbf u}(k) \subset {\mathbf U}(k) \label{eq19}\\
	Other~~constraints\label{eq20}\\
	for~~all~~k=N_{0}, N_{0}+1, ..., N \nonumber 
	\end{gather}
	Where $N$ is denoted as the last time instant in discretized version optimziation, $g_N(\mathbf x_{N})$ is the cost of final state $x_{N}$. To implement the constraints (\ref{eq18}) on state variables and other constraints (\ref{eq20}), we add penalties $\phi_N(\mathbf x_{N})$ and $\phi_k(\mathbf x_{k})$ on the cost function respectively as follows,
	
		\begin{gather}
	\min_{\substack{{\mathbf u(k)}}}~~~g_N(\mathbf x_{N})+\phi_N(\mathbf x_{N})\nonumber \nonumber \\+\sum_{N_{0}}^{N-1}(g_{k}(\mathbf x_{k},\mathbf u_{k})+\phi_k(\mathbf x_{k}, u_{k}))\times T_{s} \label{eq21}\\	
	{\mathbf x}(k+1)=\mathbf f_{k}({\mathbf x}(k),{\mathbf u}(k))\label{eq22}\\
	%{\mathbf x}(k) \subset {\mathbf X}(k) \label{eq21}\\
	{\mathbf u}(k) \subset {\mathbf U}(k) \label{eq23}\\
	%Other~~constraints\\
	for~~all~~k=N_{0}, N_{0}+1, ..., N \nonumber
	\end{gather}
	
	Finally, we bring the control objectives (\ref{eq10}) into (\ref{eq21}) and obtain the final discretized optimal control formulation.
	
	\begin{gather}
	\min_{\substack{{\mathbf u(k)}}}~~~g_N(\mathbf x_{N})+\phi_N(\mathbf x_{N})\nonumber \\ +(\sum_{{s}\in S}\sum_{N_{0}}^{N-1}a_{s}(M_{e,s}(k)-M_{e,i,d})^{2}\nonumber \\+\sum_{{i}\in N}\sum_{N_{0}}^{N-1}b_{i}|\omega_{i}(k)|\nonumber
	+\sum_{{i}\in N}\sum_{N_{0}}^{N-1}c_{i}\delta_{i}(k)^{2}\nonumber\\+\sum_{N_{0}}^{N-1}\phi_k(\mathbf x_{k}, \mathbf u_{k}))\times Ts\\
	{\mathbf x}(k+1)=\mathbf f_{k}({\mathbf x}(k),{\mathbf u}(k))\\
	%{\mathbf x}(k) \subset {\mathbf X}(k) \\
	{\mathbf u}(k) \subset {\mathbf U}(k) \\
	%Other~~constraints\\
	for~~all~~k=N_{0}, N_{0}+1, ..., N \nonumber
	\end{gather}
	
	where the first two terms in control objective are stage cost and penalty respectively for state varaibles in the final stage, and the next three terms are stage cost for state variables from stage $N_{0}$ to $N-1$, and the last term is penalties to make state variables within the constraints (\ref{eq18}) and (\ref{eq20}).
\end{itemize}

Step 3)\quad{\itshape{Run the Simulation}}
\begin{itemize}
	\item[] 
	\setlength\leftskip{4em}
	\quad After running the simulation, we obtain the optimal control input trajectory and the trajectories of state variables. Comparisons are made between cases with and without constraints, and constant inertia supply and time-variant inertia supply from storage.
\end{itemize}

\subsection{PROPT MATLAB OPTIMAL CONTROL SOFTWARE}

PROPT MATLAB OPTIMAL CONTROL SOFTWARE is a platform to solve the optimal control problem (described by ODE and DAE ) in engineering practice. Many types of control problems such as disturbance control and flight path tracking can be solved efficiently. For the detailed the information about PROPT, readers can refer to reference~\cite{rutquist2010propt}, and here we list the key steps and parameters that we need for the optimal control problem in this work.

Step 1)\quad{\itshape{Initial Parameter Setting}}
\begin{itemize}
	\item[]
	\setlength\leftskip{4em}
	\quad Determine the time interval and the length of time steps for simulation.  And then, the lower and upper limits for both the state variables and input variables should be determined. Finally, the initial and final value of state variables should be determined.
\end{itemize}

\vskip 0.1cm

Step 2)\quad{\itshape{Determine the Control Objective}}

\begin{itemize}
	\item[]
	\setlength\leftskip{4em}
	\quad To determine the control objective, both the original control objective and the penalty related to the contraints of state variables should be considered.
	
	\begin{gather}
	\min_{\substack{{\mathbf u(t)}}}\int_{t_{0}}^{t_{1}}g(\mathbf{x}(t),\mathbf{u}(t))+\phi(\mathbf{x}(t),\mathbf{u}(t))dt\label{eq27}\\
	\dot{\mathbf x}(t)=\mathbf f({\mathbf x}(t),{\mathbf u}(t))\label{eq28}\\
	%{\mathbf x}(t) \subset {\mathbf X}(t) \\
	{\mathbf u}(t) \subset {\mathbf U}(t) \label{eq29}\\
	Other~~constraints  \label{eq30}
	\end{gather}
	where $\phi(\mathbf{x}(t),\mathbf{u}(t))$ is the penalties related to the constraints of state variables (see (\ref{eq13})). 
	
\end{itemize}

Step 3)\quad{\itshape{Run the Simulation}}
\begin{itemize}
	\item[] 
	\setlength\leftskip{4em}
	\quad After running the simulation, we obtain the optimal control input trajectory and the trajectories of state variables. 
\end{itemize}
\vskip 0.3cm

%Step 4)\quad{\itshape{Termination}}
%\begin{itemize}
%	\item[] 
%	\setlength\leftskip{4em}
%	\quad When the on-going algorithm meets the termination condition such as the maximum numbers of generation, the calcualtion stops and returns the current best fitness value and the corresponding individual vector (decision variables defined in Step 1).
%\end{itemize}

%For calculation of $\omega_{p}^{max}$ and $\dot{\omega_{p}}^{max}$ in objective function, inspired by the reference~\cite{zhu2018optimization}, we adopt a simulation-based approach and sample the curve determined by (\ref{eq18}), and the sampling time is 0.1s. This method can avoid the problem that the obtained value of $\omega_{p}^{max}$ and $\dot{\omega_{p}}^{max}$ will just be a local optimal point by the sensitivity method in~\cite{borsche2015effects, borsche2017placement}. For the coefficients $c_{1}$-$c_{8}$ in the objective function,  their values are equal to 1 if there is no special statement.

%For calculation of $\omega_{p}^{max}$ and $\dot{\omega_{p}}^{max}$ in objective function, we sample the curve determined by (\ref{eq18}), and the sampling time is 0.1s. This method can avoid the problem that the obtained value of $\omega_{p}^{max}$ and $\dot{\omega_{p}}^{max}$ will just be a local optimal point by the sensitivity method in~\cite{borsche2015effects, borsche2017placement}. For the coefficients $c_{1}$-$c_{8}$ in the objective function,  their values are equal to 1 if there is no special statement.

\section{Case study}

We will do two case studies in this section. The first case study is conducted in a simple two-bus system, where typical grid parameters are utilized. And next, we will employ a more realistic 12-bus system to verify the solving method. 

For the case study in 2-bus system, the simulation is run in Matlab R2018a on a computer with Microsoft Windows, %Version 1607, i5 6500 3.2 GHz processor, 1TB hard disk and 16GB RAM. 
i5-7500 3.40GHz, and 8GB RAM. For the Matlab code to implement the basic DP and level-set DP, readers can refer to reference \cite {sundstrom2009generic} and download it at \cite {kkk}. For parameter initialization, it is suggested by authors that the discretization of state variable space and time step should be as fine as possible if computation ability is allowed. Or there is a compromise between your discretizaiton step and your computation ability. The result for comparison is obtained from LS DP, if there is no special statement.

For the case study in 12-bus system, the PROPT MATLAB OPTIMAL CONTROL SOFTWARE of trial version is utilized, and it can be downloaded at~\cite{Sdd}. And `High Dim Control' is selected for `options.name' in the optimization program. The simulation for this case is run on MacBook Pro, 2.3GHz Intel Core i5 and 8GB 2133MHz LPDDR3 memory. 

\subsection{2-bus system}

\begin{figure}[!tbh]
	\centering
	\includegraphics[width=2in]{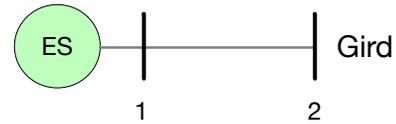}
	\caption{2-bus test system.}
	\label{n1}
\end{figure}

Fig.~\ref{n1} shows the a two-bus system where storage at bus 1 connects to the grid at bus 2. Bus 2 is set as a reference bus which means $\delta_{2} \equiv0$ and $\omega_{2}\equiv0$ during the transient process. The initial parameters are given in Tab.~\ref{tn11} for simulation, we can see that the setting is that there is no power exchange between the storage and the grid at the initial time instant, and the disturbance is set as the  0.3 p.u. power increase at bus $1$.

\begin{table}[!tbh]
	\renewcommand{\arraystretch}{1.}
	\caption{Initial parameters for simulation for 2-bus system.}
	\label{tn11}
	\centering
	\begin{tabular}{cc|cc}
		\hline
		\hline
		\textbf{Parameter}  &  \textbf{Value}& \textbf{Parameter}&\textbf{Value}\\
		\hline
		$T_{s}$  &  0.5 s& $t_{1}$  &  30 s\\ 
		\hline
		$N_{\delta_{1}}$ & 201  &$N_{\omega_{1}}$ &   51\\
		\hline
		$\delta_{1}(0)$ & 0 & $\omega_{1}(0)$ & 0 \\
		\hline
		$\delta_{1}(t)$ & [0, 0.6]& $\omega_{1}(t)$ & [-0.5, 0.5]\\
		\hline
		$\delta_{1}(N)$ & [0, 0.6] & $\omega_{1}(N)$ & [-0.02, 0.02]\\
		\hline
		$M_{1}(t)$ & [4 s, 10 s]& $B_{1-2}$ & 1 \\
		\hline
		$D_{e,1}$  & 1  & $N_{M_{e,1}}$  &  51\\
		\hline
		\hline
	\end{tabular}
\end{table}

\subsubsection{Time-variant Inertia vs. Constant Inertia}
\begin{figure}[!tbh]
	\centering
	\includegraphics[width=3.35in]{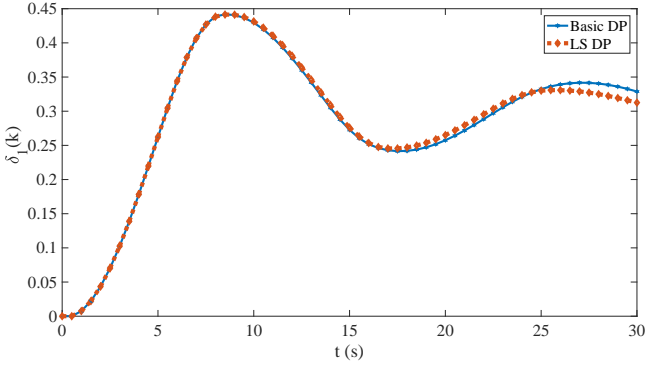}
	\caption{Angle at bus 1.}
	\label{m1}
\end{figure}

\begin{figure}[!tbh]
	\centering
	\includegraphics[width=3.35in]{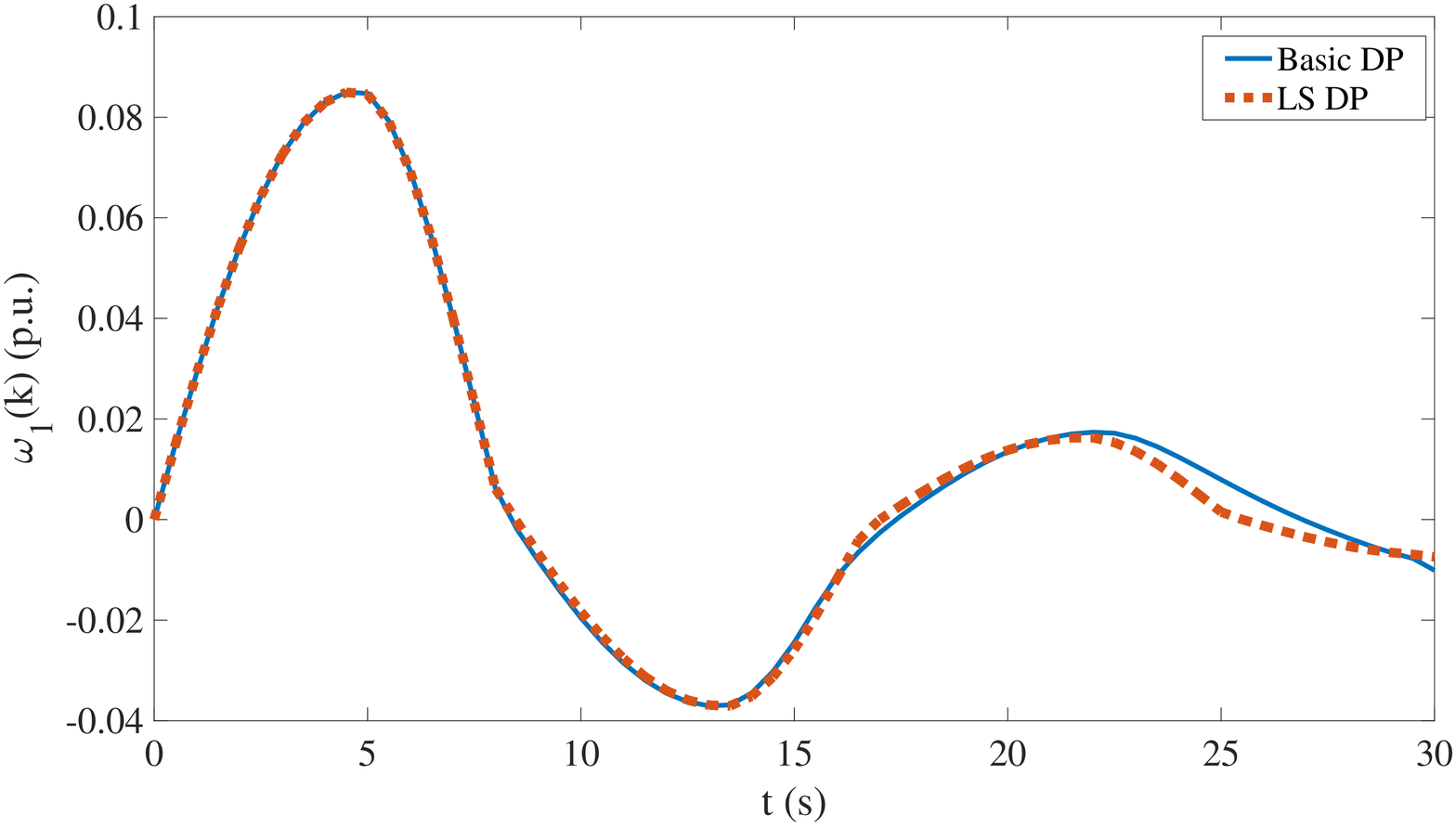}
	\caption{Frequency at bus 1.}
	\label{m2}
\end{figure}

\begin{figure}[!tbh]
	\centering
	\includegraphics[width=3.35in]{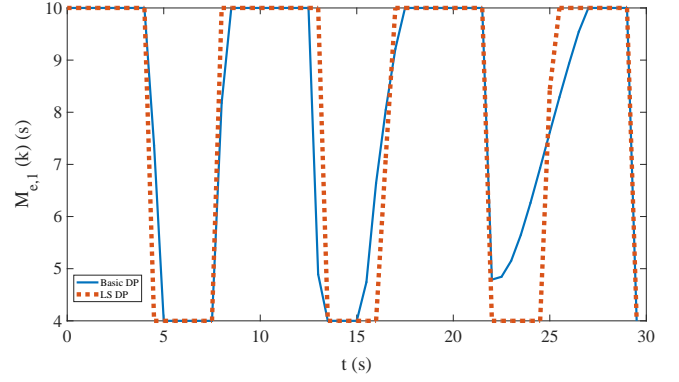}
	\caption{Inertia change at bus 1.}
	\label{m3}
\end{figure}

First we set the control objective is as follows:
\begin{gather}
	\min_{\substack{{\mathbf u(k)}}}~~~\phi_N(\mathbf x_{N})+\sum_{{i}\in N}\sum_{N_{0}}^{N-1}b_{i}|\omega_{i}(k)|
\times Ts\label{eq31}
\end{gather}
where $b_{i}$ is the coefficient and equals to 1, and $\phi_N(\mathbf x_{N})$ is 2 here. After running the simulation, we obtain the optimal trajectories of $\delta_{1}(k)$, $\omega_{1}(k)$ and $M_{e,1}(k)$, as shown in Fig~\ref{m1}, Fig~\ref{m2} and Fig~\ref{m3} respectively. %For the comparison in following case studys, the value is based on the LS DP algorithm if no special illustation is made.
\vskip 0.3cm

For the basic DP, the value of objective function is 2.7549 (frequency absolute value integration: 0.7549 + penalty: 2), and for the LS DP, the value is 2.7408 (frequency absolute value integration: 0.7508 + penalty: 2), it can be seen that LS DP can achieve better result compared with basic DP. 

We also see some interesting phenomena when the inertia can be changed with time:

   \begin{enumerate}[i.]
	\item When the time $t$ is between 0s and 5s, the maximum virtual inertia is chosen by the storage, this is to prevent the increase of frequency $\omega_{1}$.
	\item When the value of frequency $\omega_{1}$ changes sign, or the frequency changes the direction of motion, the inertia change will switch between two boundary values, this is also to prevent the frequency change.
	\item The value of frequency at the final state is within the predefined range, which means penalty $\phi_{N}(\mathbf x_{N})$ takes effects.
\end{enumerate}

The results listed above show that this optimal control problem has been successfully solved. To compare this case with one where inertia is fixed, we adopt the following control objective,

\begin{gather}
\min_{\substack{{\mathbf u(k)}}}~~~\phi_N(\mathbf x_{N})+\sum_{{i}\in N}\sum_{N_{0}}^{N-1}b_{i}|\omega_{i}(k)|
\times Ts \nonumber \\
+\sum_{{s}\in S}\sum_{N_{0}}^{N-1}a_{s}(M_{e,s}(k)-M_{e,i,d})^{2}\times T_{s} \label{eq32}
\end{gather}
where $a_{s}$ equal 100000, and desired inertia $M_{e,i,d}$ is 4s. Through this control objective, we want to see that the virtual inertia will be at 4s over the time.

\begin{figure}[!tbh]
\centering
\includegraphics[width=3.35in]{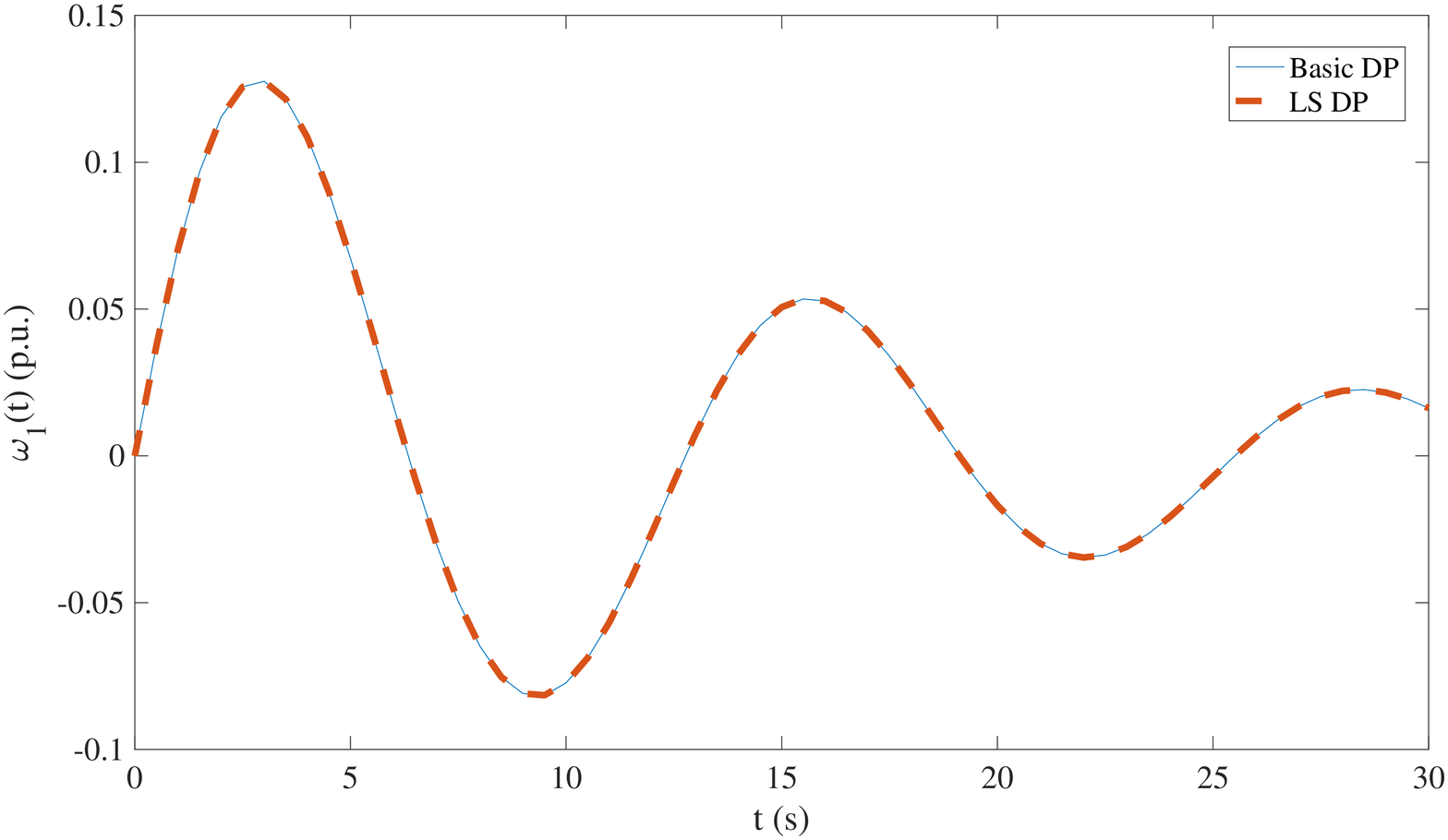}
\caption{Frequency at bus 1.}
\label{mn1}
\end{figure}

\begin{figure}[!tbh]
\centering
\includegraphics[width=3.35in]{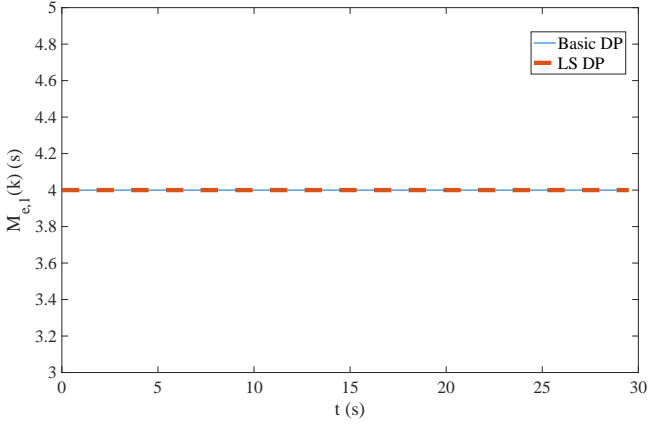}
\caption{Inertia change at bus 1.}
\label{m5}
\end{figure}

It can seen from Fig.~{\ref{mn1}} and Fig.~{\ref{m5}} that inertia at bus 1 keeps constant at 4s and frequency at bus 1 experiences larger oscillations than that in Fig.~{\ref{m2}}. For the frequency absolute value integration with time in contant-inertia case, the value is 1.2792. The running time of basic DP and LS DP for control objective (\ref{eq31}) is 91.6709s and 115.4804s, and the running time of basic DP and LS  DP for control objective (\ref{eq32}) is 84.4507s and 97.7915s respectively.

\subsubsection{Verifying the other constraints} In this subsection, the author would like to verify that dynamic programming can meet the other constraints such as the energy constraints of storage. The case (control objective (\ref{eq32})) in the last subsection is a base case. We want to add the following constraint: 

\begin{gather}
P_{i}^{r,min}\leq P_{i}^{r}= P_{i}^{e}-M_{e,i}\dot{w^{e}}_{i}-D_{e,i}w^{e}_{i}\leq P_{i}^{r,max}\label{eq33}\\
E^{al,l}_{i}\leq \int_{t_{0}}^{t_{1}}P_{i}^{r} dt\leq E^{al,u}_{i}\label{eq34}
\end{gather}

They are power capacity constraint and energy capcacity cosntraints respectively. And we bring the parameters into these two constraints and discretize them as follows,

\begin{gather}
P_{i}^{r,min}\leq P_{i}^{r}(k)=- M_{e,i}\frac{{w^{e}}_{i}(k+1)-{w^{e}}_{i}(k)}{T_{s}}\nonumber \\-w^{e}_{i}(k)\leq P_{i}^{r,max} \label{eq35}\\
E^{al,l}_{i}\leq \sum_{N_{0}}^{N-1}P_{i}^{r} (k)\leq E^{al,u}_{i}\label{eq36}
\end{gather}

\begin{figure}[!tbh]
	\centering
	\includegraphics[width=3.35in]{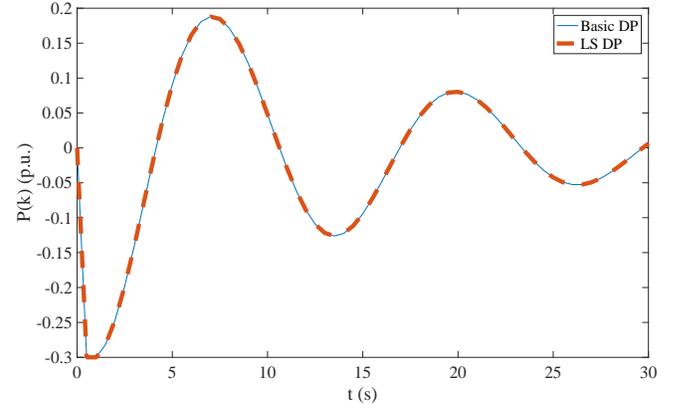}
	\caption{Power change of storage at bus $1$.}
	\label{m7}
\end{figure}

\begin{figure}[!tbh]
	\centering
	\includegraphics[width=3.35in]{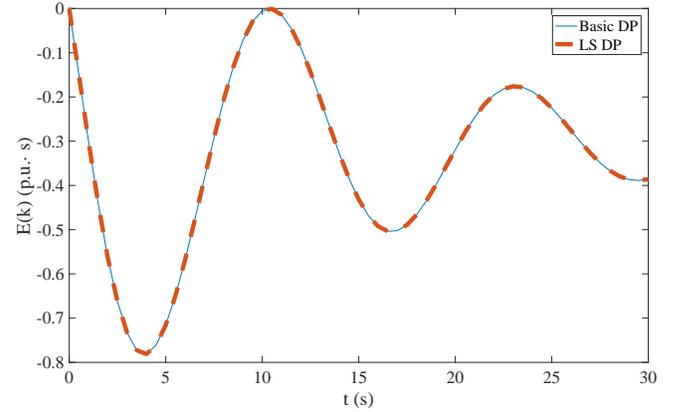}
	\caption{Energy change of storage at bus 1.}
	\label{m8}
\end{figure}

\begin{figure}[!tbh]
	\centering
	\includegraphics[width=3.35in]{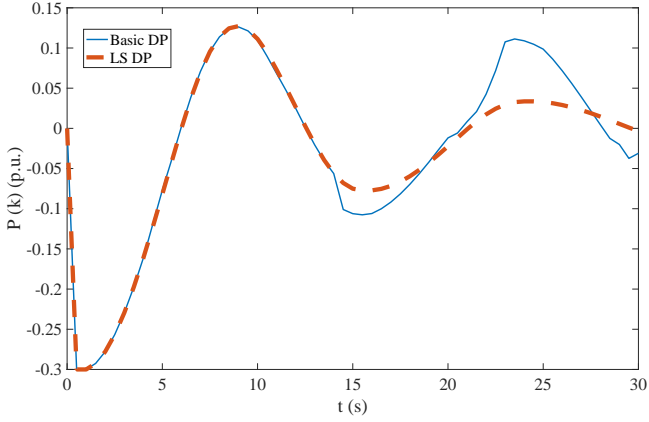}
	\caption{Power change of storage at bus 1.}
	\label{m9}
\end{figure}

\begin{figure}[!tbh]
	\centering
	\includegraphics[width=3.35in]{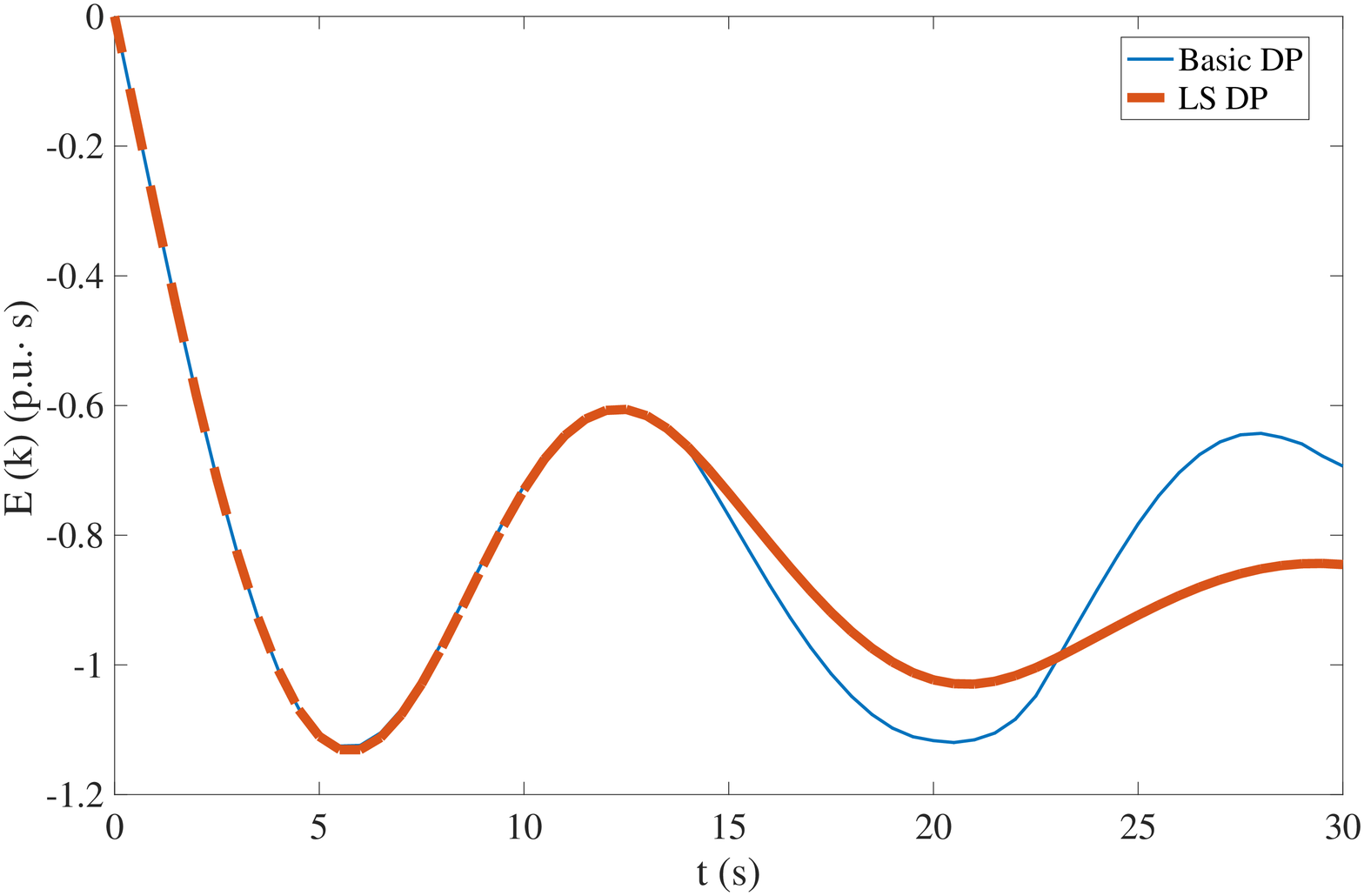}
	\caption{Energy change of storage at bus 1.}
	\label{m10}
\end{figure}

\begin{figure}[!tbh]
	\centering
	\includegraphics[width=3.35in]{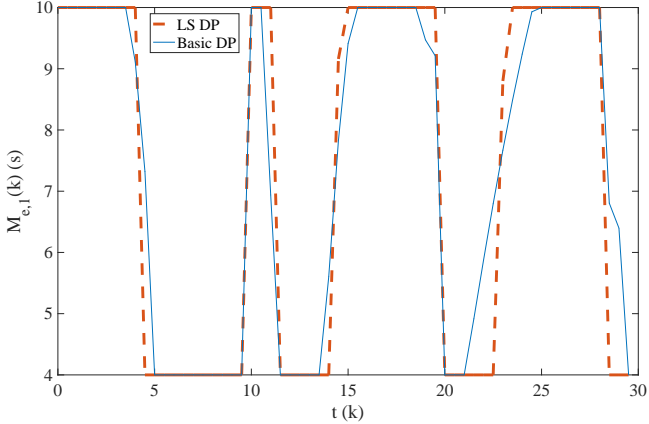}
	\caption{Virtual inertia change at bus 1.}
	\label{m11}
\end{figure}

For the base case, the power change and energy change of storage are shown in Fig.~\ref{m7} and in Fig.~\ref{m8} respectively. The minimum and maximum value of $P(k)$ is -0.3 p.u and 0.1880 p.u. at $t$=0.5s and at $t$=7s respectively. It is noted that we utilize the  discretization method and intial value of state varables are 0, so the minimum value of $P(k)$ at the first step $t$=0.5s is fixed and is equal to -0.3 p.u.. Another fact is that this is a single-machine infinite-bus system, the controllability of state variables are limited.  And thus, we set $P_{i}^{r,max}=0.15$ p.u. And we penalize the power change by upper limits, this will significantly increase the energy in the storage. To make the result more clear, we only implement power constraints for storage, and the control objective is as follows, % The minimum value of $E(k)$ is -0.7812 p.u.$\cdot$s at time  $t$=4s and we set $E_{i}^{r,min}=-0.7$ p.u.$\cdot$s.
\begin{gather}
\min_{\substack{{\mathbf u(k)}}}~~~\phi_N(\mathbf x_{N})+\sum_{{i}\in N}\sum_{N_{0}}^{N-1}b_{i}|\omega_{i}(k)|
\times Ts \nonumber \\
+\sum_{{s}\in S}\sum_{N_{0}}^{N-1}d_{s}max(P^{r}_{i}(k)-0.15, 0)\label{eq37}
\end{gather}
where $d_{s}$ is the coefficient, and equals to 100000, and $b_{i}$ equals to 1. The power change, energy change of storage, and virtual inertia trajectory are shown in Fig.~\ref{m9}, Fig.~\ref{m10} and Fig.~\ref{m11} respectively. It can be seen that power change of storage does not go beyond its upper limit 0.15 p.u. during the transient process and the maximum value of $P(k)$  is 0.1272 p.u. at time $t$=9s. And running time for basic DP and LS DP is 82.7141s and 90.8032s respectively.

In this optimal control problem, constraint (\ref{eq36}) is related to multi-stage state variables, in fact, it is noted that storage capacity can also be treated as a state variable, and this constraint will be easy to be transformed into a penalty function, which can be added on the control objective. Since this paper focuses on the effect of time-variant virtual inertia, we only point out the feasible methods to deal with these kinds of constraints.

\subsection{12-bus system}
\begin{figure}[!tbh]
	\centering
	\includegraphics[width=3.35in]{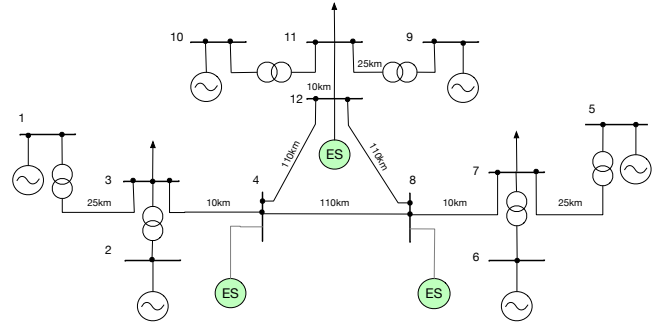}
	\caption{12-bus test system.}
	\label{f2}
\end{figure}

The 12-bus test system in Fig.~\ref{f2} is modified from the well-known two-area system in reference~\cite{kundur1994power} and an additional area is added as reference~\cite{borsche2015effects}.% It contains 6 generators and 6 loads. 
~The transformer reactance is 0.15 p.u. and the line impedance is (0.0001+0.001i) p.u./km. We still utilize structure preserving model to describe the dynamics of power system. The base capacity of this system for power flow calculation is set as 100MVA. The inertia and damping of original power system is given in Table~\ref{tab1} and the steady power flow condition is given in Table~\ref{tab22}.  It is assumed that there are motor loads (including little inertia and damping) at the load buses. And bus 9 is a set as a reference bus in the system. The time step is 0.5s and the time interval for running the simulation is 40s. Our eyesight in this case will be put on minimizing frequency devations and power flow oscillations on transmission lines. The contingency setting is power increase of 60MW (0.6 p.u.) at bus 1.

%\begin{table}[!tbh]
%	\renewcommand{\arraystretch}{1.}
%	\caption{Initial Power Flow Condition.}
%	\label{tab2}
%	\centering
%	\begin{tabular}{c|ccccccc}
%		\hline
%		\hline
%		Gen  &  1  &  2  &  5 & 6  & 9 & 10\\
%		\hline
%		$P$ (MW)  & 138& 1050 & 719 & 350 & 700 & 700\\
%		\hline
%		Load  &  3  & 4  &  7 & 8  & 11 & 12\\
%		\hline
%		$P$ (MW)  &  400& 567 & 490 & 800 & 400 & 1000\\
%		\hline
%		\hline
%	\end{tabular}
%\end{table}

%\begin{table}[!tbh]
%	\renewcommand{\arraystretch}{1.}
%	\caption{Inertia and Damping Distribution of Original Power System.}
%	\label{tab1}
%	\centering
%	\begin{tabular}{cc}
%		\hline
%		\hline
%		\textbf{Bus. No.}  &  \textbf{Inertia (s) / Damping (p.u.)}\\
%		\hline
%		1, 2  &  15/3 \\
%		\hline
%		3, 7, 11 & 1/0.1\\
%		\hline
%		5, 6 & 20/4\\
%		\hline
%		9, 10& 10/2\\
%		\hline
%		\hline
%	\end{tabular}
%\end{table}
\subsubsection{Base Case} 
\begin{figure}[!tbh]
	\centering
	\includegraphics[width=3.35in]{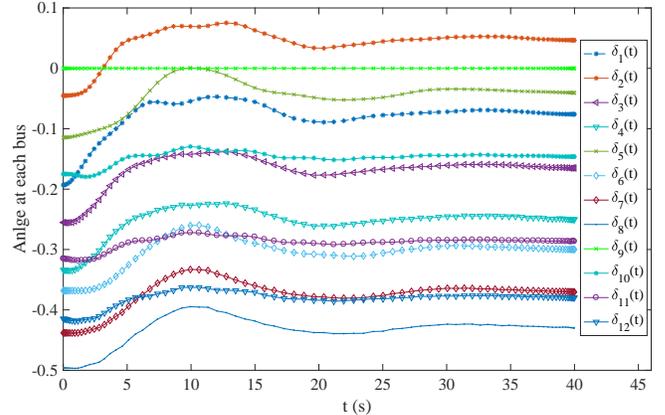}
	\caption{Angle at each bus.}
	\label{mm4}
\end{figure}

\begin{figure}[!tbh]
	\centering
	\includegraphics[width=3.35in]{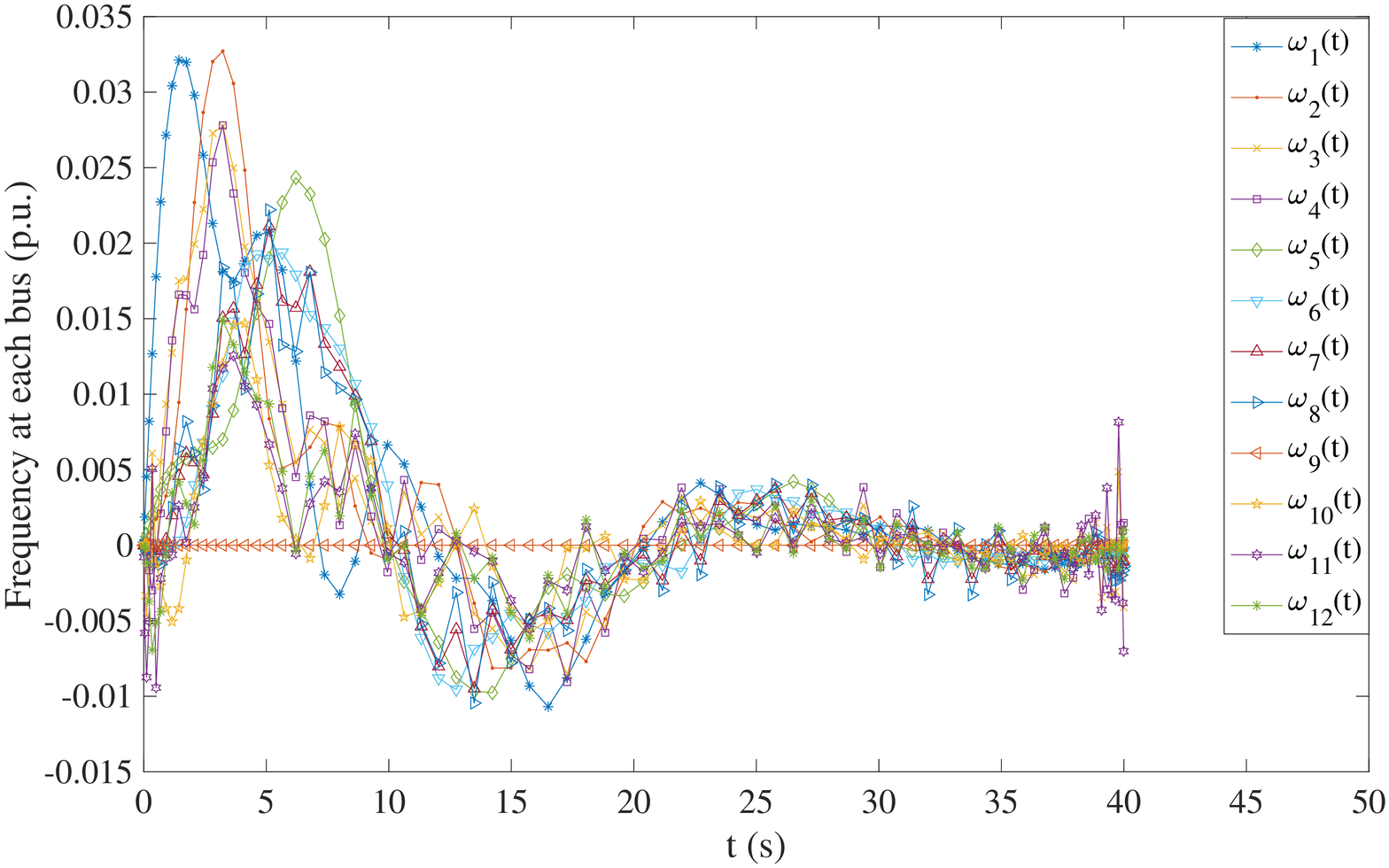}
	\caption{Frequency at each bus.}
	\label{mm5}
\end{figure}

\begin{figure}[!tbh]
	\centering
	\includegraphics[width=3.35in]{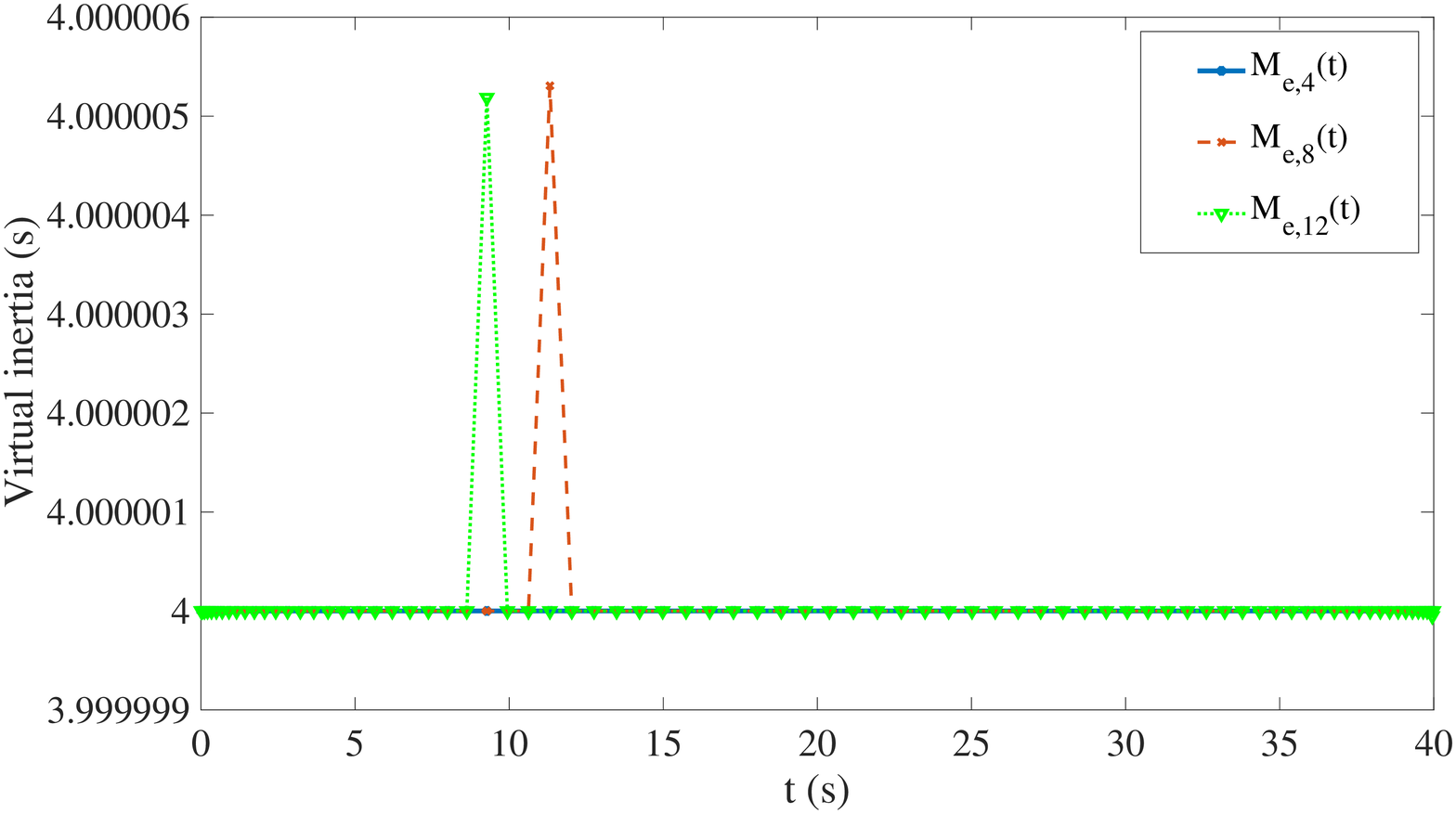}
	\caption{Virtual inertia change at each bus.}
	\label{mm6}
\end{figure}

\begin{figure}[!tbh]
	\centering
	\includegraphics[width=3.35in]{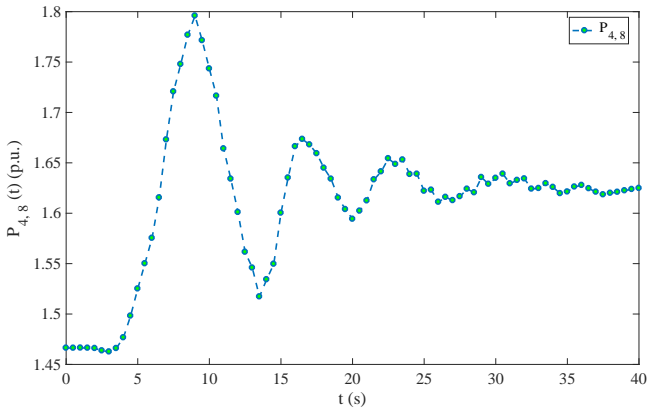}
	\caption{Power flow at transmission line 4-8.}
	\label{mmf}
\end{figure}

For comparison purpose, we first do a base case where the virtual inertia of storage is fixed. And we adopt $a_{s}$=1 and $M_{e,i,d}=4$ for $i\in S$, $b_{i}$=0, $c_{i}$=0
for $i\in N$ in control objective (\ref{eq10}). And the angle, frequency,  inertia trajectories and power flow from bus 4 to bus 8 are shown in Fig.~\ref{mm4}, Fig.~\ref{mm5}, Fig.~\ref{mm6} and Fig.~\ref{mmf}.

The time for this case study is 1.8123s for symbolic processing and 3.9842s for CPU calculation. The time integration for frequency absolute value is 1.8157 p.u.$\cdot$s. The power peak is 1.7959 p.u. at $t$=9s.

\subsubsection{Minimizing frequency deviations}

\begin{figure}[!tbh]
	\centering
	\includegraphics[width=3.35in]{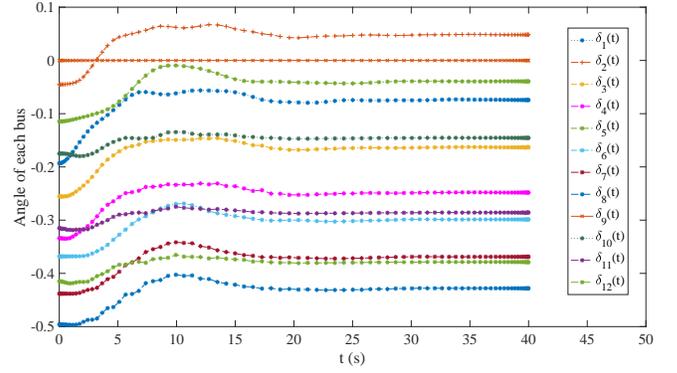}
	\caption{Angle at each bus.}
	\label{mm1}
\end{figure}

\begin{figure}[!tbh]
	\centering
	\includegraphics[width=3.35in]{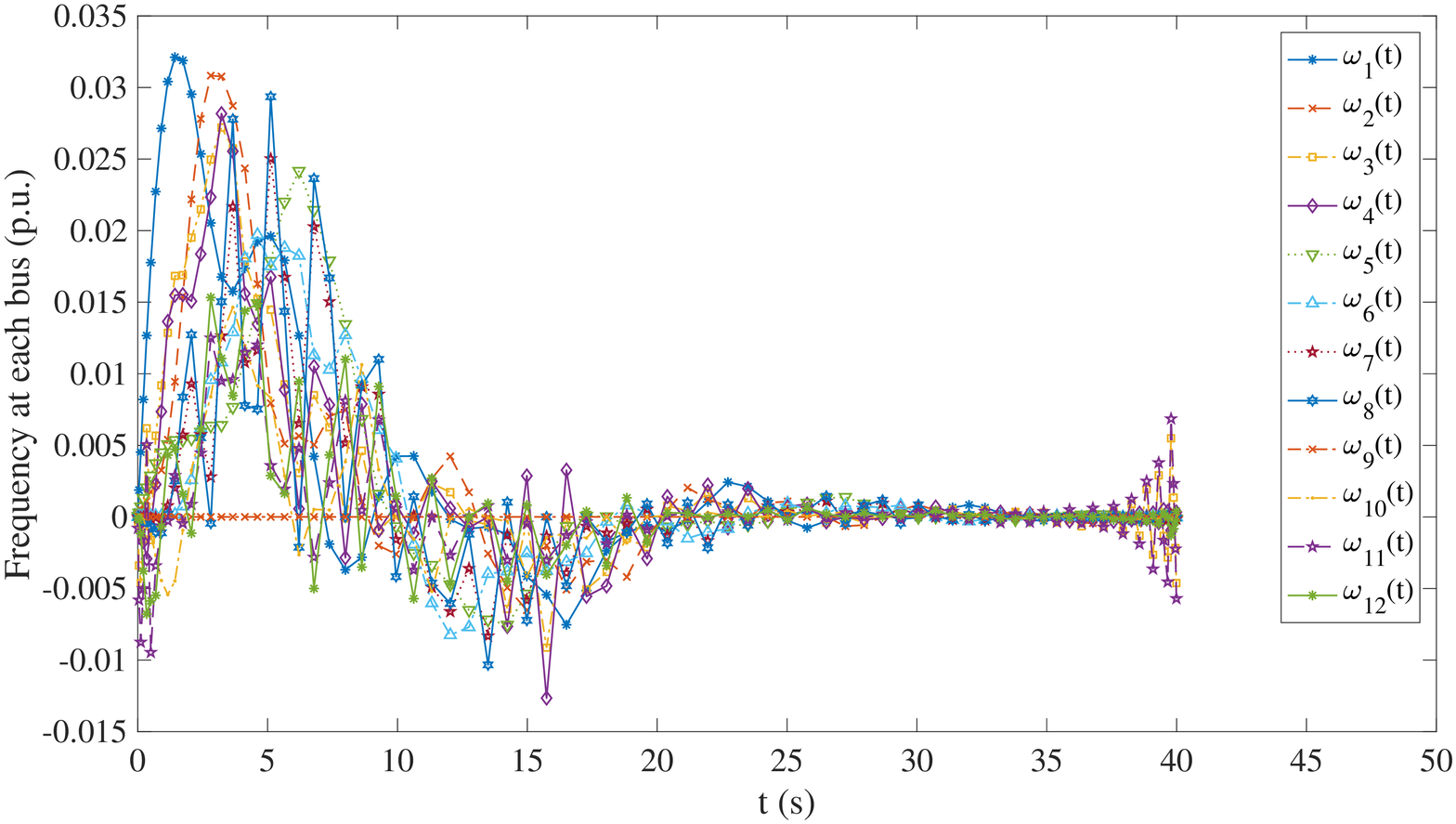}
	\caption{Frequency at each bus.}
	\label{mm2}
\end{figure}

\begin{figure}[!tbh]
	\centering
	\includegraphics[width=3.35in]{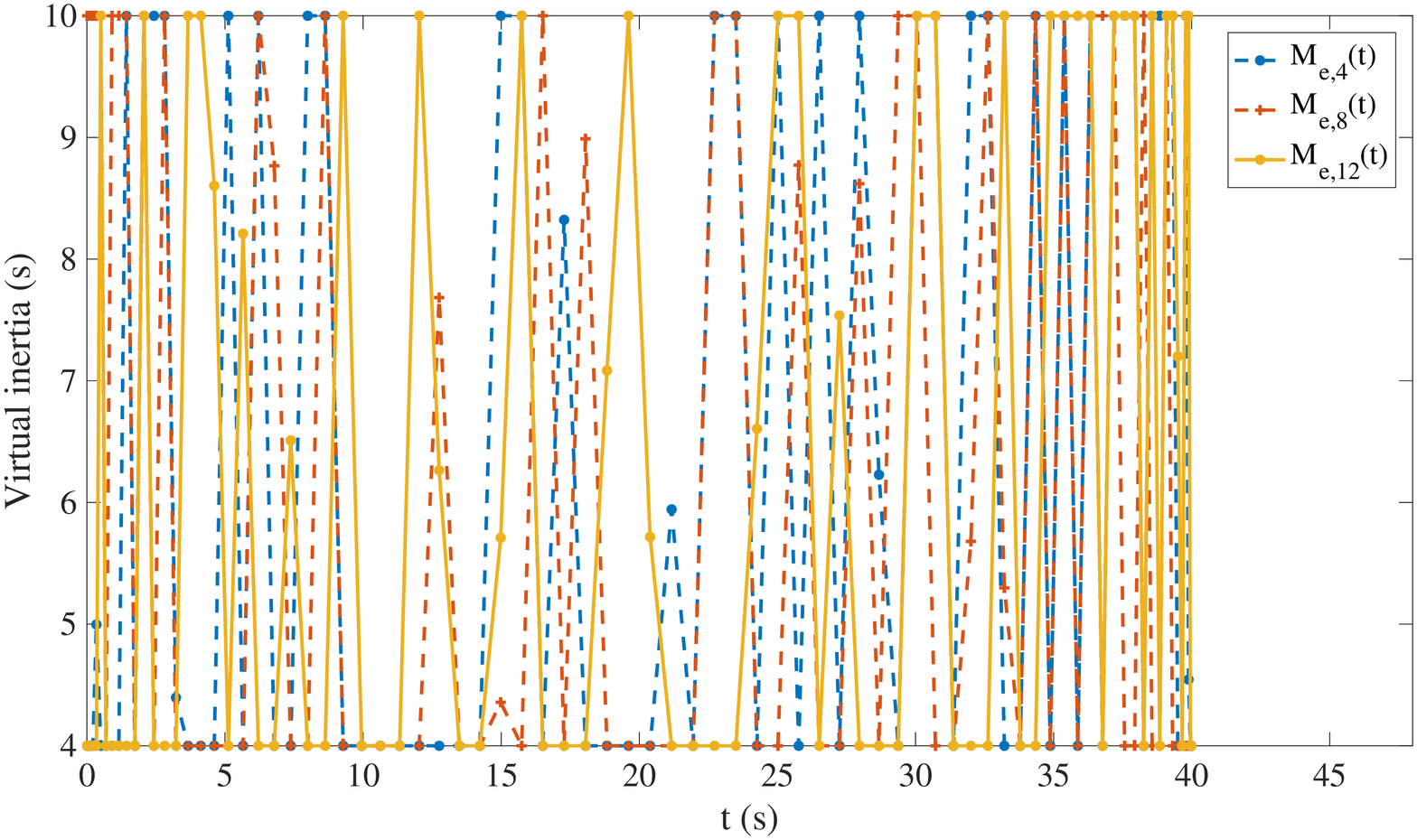}
	\caption{Virtual inertia change at each bus.}
	\label{mm3}
\end{figure}

To minimize the frequency deviations, we do a case where $a_{s}$=0 for $i\in S$ , $b_{i}$=1, $c_{i}$=0
for $i\in N$. And the angle, frequency and inertia trajectories are shown in Fig.~\ref{mm1}, Fig.~\ref{mm2} and Fig.~\ref{mm3} as follows, 

The time for the base case study is 3.5128s for symbolic processing and 247.65s for CPU calculation. The time integration for frequency absolute value is 1.4975 p.u.$\cdot$s, and we can see that the control objective for frequency minimization is achieved.

\subsubsection{Minimizing the power flow oscillations} To minimize the power flow ocillations, we adopt the following control objective,
\begin{gather}
\label{eq24}
\min_{\substack{{M_{e,s}(t)}}}\int_{t_{0}}^{t_{1}}max(b_{48}sin(\delta_{4}-\delta_{8})-1.7, 0)~dt
\end{gather}

\begin{figure}[!tbh]
 	\centering
 	\includegraphics[width=3.35in]{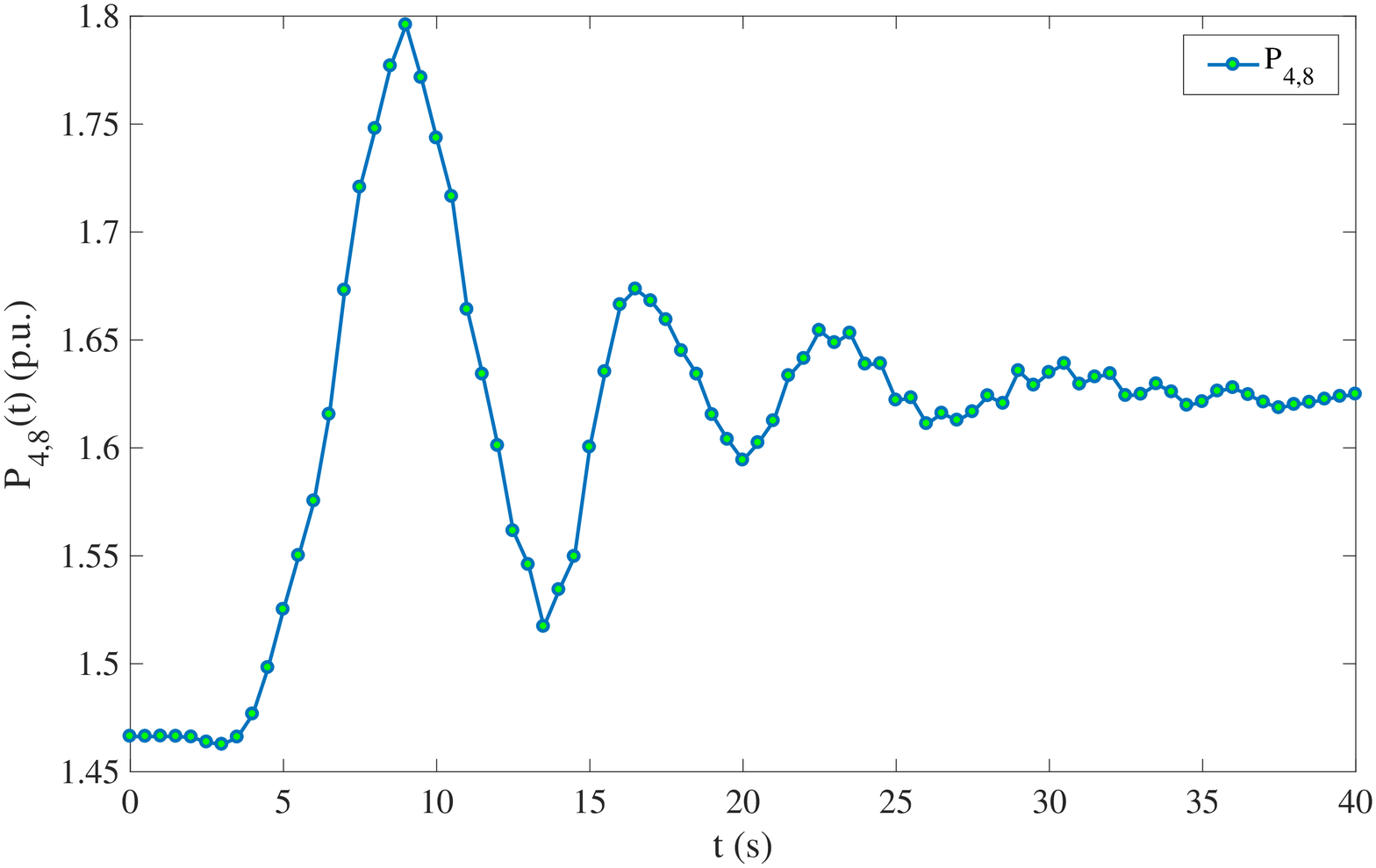}
 	\caption{Power flow at transmission line 4-8.}
 	\label{mmf1}
 \end{figure}

\begin{figure}[!tbh]
	\centering
	\includegraphics[width=3.35in]{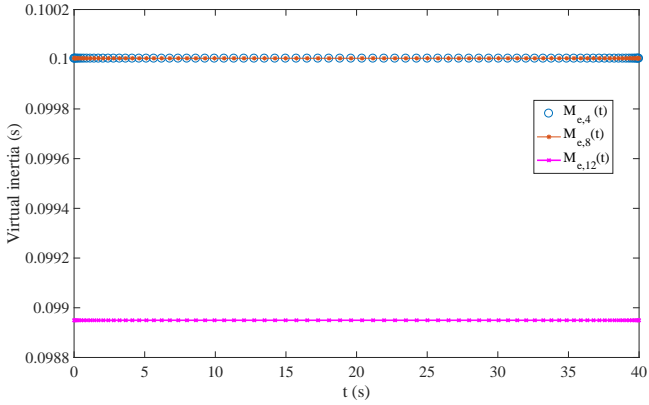}
	\caption{Virtual inertia change at each bus.}
	\label{mmf2}
\end{figure}

 We enlarge the virtual inertia range to [0.1s, 15s], and the power flow and virtual inertia trajectories are shown in Fig.~\ref{mmf1} and Fig.~\ref{mmf2}.
 
The time for this case study is 3.0355s for symbolic processing and 4.3600s for CPU calculation.  The power peak is 1.7853 p.u. at $t$=9s. And we can see that even though we enlarge the virtual inertia range, in this case, the effect of virtual inertia change on the minimizing the power flow oscillation is limited. The virtual inertia provided by all storage is the minimum value at around 0.1s.

\section{conclusion}

This paper novelly treats controlling time-variant virtual inertia as an optimal control problem, and provides two corresponding methods, dynamic programming and PROPT respectively, to solve it. Dynamic programming is a generally applicable method for different kinds of system models, which means the analysis can be conducted for power system with different types of power-electronic equipments. For the PROPT method, this software can deal with high-dimensional control with very fast speed, which is desired for analysis.

This work opens a new space for frequency control,  research opportunities followed by this work can be as follows,

   \begin{enumerate}[i.]
	\item Dynamic programming will suffer the problem `curse of dimensionality', the calculation time will exponentially increase with the increase of state variables. For the high dimensional dynamic system, what should we do if we still want to use dynamic programming method?
	\item In this work, the wind power change is not considered, if wind power is considered, how do we change the model correspondingly?
	\item In this work, we only treat the virtual inertia as the control input, in fact, power input/output of storage can also be treated as the control input with a certain cost coefficient, so which is better between controlling power or controlling virtual inertia to achieve a specific control objective?
\end{enumerate}

\newpage

%This paper analyzes the problem that where and how much virtual inertia and damping is emulated for energy storage in power system by adopting a model-based method. Differential evolution algorithm is employed to solve this mix-integer nonlinear problem. The results verify the effectiveness of this algorithm, and damping ratio, frequency overshoot, RoCoF and energy change of energy storage are optimized respectively in case study. For the former three objectives, the optimal solutions are not on the boundary of domain of decision variables, and this validates the nonlinearities of the problem. 
%
%Due to the fact that linearized model in this work is only effective for small signal stability analysis. For the futher work, authors would like to analyze the question that where and how much virtual inertia and damping is deployed for energy storage can enhance the stability of power system when a large disturbance happens. 

\vskip 1cm

% or


\appendices
\section{Initial Parameters for 12-bus system simulation.}

\begin{table}[!tbh]
	\renewcommand{\arraystretch}{1.}
	\caption{Power Flow Condition for 12-bus system.}
	\label{tab22}
	\centering
	\begin{tabular}{c|ccccccc}
		\hline
		\hline
		Gen  &  1  &  2  &  5 & 6  & 9 & 10\\
		\hline
		$P$ (MW)  & 138& 1050 & 719 & 350 & 700 & 700\\
		\hline
		Load  &  3  & 4  &  7 & 8  & 11 & 12\\
		\hline
		$P$ (MW)  &  400& 567 & 490 & 800 & 400 & 1000\\
		\hline
		\hline
	\end{tabular}
\end{table}

\begin{table}[!tbh]
	\scriptsize
	\renewcommand{\arraystretch}{1.}
	\caption{Initial parameters for simulation of 12-bus system.}
	\label{tn1}
	\centering
	\begin{tabular}{cc|cc}
		\hline
		\hline
		\textbf{Parameter}  &  \textbf{Value}& \textbf{Parameter}&\textbf{Value}\\
		\hline

		\hline
		$\delta_{1}(0)$ & -0.1931 & $\omega_{1}(0)$ & 0 \\
		\hline
		$\delta_{2}(0)$ & -0.0452 & $\omega_{2}(0)$ & 0 \\
		\hline
		$\delta_{3}(0)$ & -0.2552 &$\omega_{3}(0)$ & 0 \\
		\hline
		$\delta_{4}(0)$ & -0.3340 & $\omega_{4}(0)$ & 0 \\
		\hline
		$\delta_{5}(0)$ & -0.1146 & $\omega_{5}(0)$ & 0 \\
		\hline
		$\delta_{6}(0)$ & -0.3681 & $\omega_{6}(0)$ & 0 \\
		\hline
		$\delta_{7}(0)$ & -0.4381 & $\omega_{7}(0)$ & 0 \\
		\hline
		$\delta_{8}(0)$ & -0.4960 & $\omega_{8}(0)$ & 0 \\
		\hline
		$\delta_{9}(0)$ & 0 &  $\omega_{9}(0)$  & 0 \\
		\hline
		$\delta_{10}(0)$ & -0.1750 &$\omega_{10}(0)$ & 0 \\
		\hline
		$\delta_{11}(0)$ & -0.3150 & $\omega_{11}(0)$ & 0  \\
		\hline
		$\delta_{12}(0)$ & -0.4150 & $\omega_{12}(t)$ & 0  \\
		\hline
		$D_{e,4}$  & 0.1 p.u.  & $D_{e,8}$  &  0.1 p.u. \\
			\hline
		$D_{e,4}$  & 0.1 p.u. &  $M_{e,4}(t)$ &  [4s, 10s] \\
		\hline 
		$M_{e,8}(t)$  & [4s, 10s]  &  $M_{e,12}(t)$ &  [4s, 10s] \\
		\hline
		\hline
	\end{tabular}
\end{table}

\begin{table}[!tbh]
	\renewcommand{\arraystretch}{1.}
	\caption{Inertia and Damping Distribution of Original Power System.}
	\label{tab1}
	\centering
	\begin{tabular}{cc}
		\hline
		\hline
		\textbf{Bus. No.}  &  \textbf{Inertia (s) / Damping (p.u.)}\\
		\hline
		1, 2  &  15/3 \\
		\hline
		5, 6 & 20/4\\
		\hline
		9, 10& 10/2\\
		\hline
		3, 7, 11  & 1/0.1\\
		\hline
		\hline
	\end{tabular}
\end{table}

\ifCLASSOPTIONcaptionsoff
  \newpage
\fi
\newpage



\bibliographystyle{IEEEtran}
\bibliography{journal}

\end{document}